\PassOptionsToPackage{unicode}{hyperref}
\PassOptionsToPackage{hyphens}{url}
\PassOptionsToPackage{dvipsnames,svgnames,x11names}{xcolor}
\documentclass[
12pt,letterpaper]{article}

\RequirePackage{amsthm,amsmath,amsfonts,amssymb}
\RequirePackage[authoryear,round]{natbib}
\bibpunct{(}{)}{;}{a}{}{,}
\RequirePackage[colorlinks,citecolor=blue,urlcolor=blue]{hyperref}
\usepackage{dutchcal}
\usepackage{booktabs}   
\usepackage{caption}
\usepackage{multirow}
\usepackage{bm}
\usepackage{iftex}
\newif\ifJASAUnicodeTeX
\ifXeTeX\JASAUnicodeTeXtrue\fi
\ifLuaTeX\JASAUnicodeTeXtrue\fi
\ifJASAUnicodeTeX
\usepackage{unicode-math}
\defaultfontfeatures{Scale=MatchLowercase}
\defaultfontfeatures[\rmfamily]{Ligatures=TeX,Scale=1}
\else 
\usepackage[T1]{fontenc}
\usepackage[utf8]{inputenc}
\usepackage{textcomp} 
\fi
\usepackage{lmodern}
\newif\ifJASAMicrotype
\ifPDFTeX\JASAMicrotypetrue\fi
\ifJASAUnicodeTeX\JASAMicrotypetrue\fi
\ifJASAMicrotype
\IfFileExists{microtype.sty}{
	\usepackage[]{microtype}
	\UseMicrotypeSet[protrusion]{basicmath} 
}{}
\fi
	\IfFileExists{parskip.sty}{%
		\usepackage{parskip}
	}{
		\setlength{\parindent}{0pt}
		\setlength{\parskip}{6pt plus 2pt minus 1pt}}
\usepackage{xcolor}
\usepackage{array}
\usepackage[letterpaper,margin=1in]{geometry}
\ifLuaTeX
\usepackage{selnolig}  
\fi
\usepackage{bookmark}

\IfFileExists{xurl.sty}{\usepackage{xurl}}{} 
\hypersetup{
	pdftitle={When Standardization Spreads Factor Variation Across the Spectrum},
	pdfauthor={},
	pdfkeywords={Factor variation, high dimensions, sample correlation matrices, diverging spikes},
	colorlinks=true,
	linkcolor={blue},
	filecolor={Maroon},
	citecolor={Blue},
	urlcolor={Blue},
	pdfcreator={LaTeX}}
\theoremstyle{plain}
\newtheorem{theorem}{Theorem}[section]
\newtheorem{lemma}[theorem]{Lemma}
\theoremstyle{definition}
\newtheorem{definition}[theorem]{Definition}
\newtheorem{example}[theorem]{Example}
\newtheorem{assumption}{Assumption}
\newtheorem{proposition}[theorem]{Proposition} 
\newtheorem{corollary}[theorem]{Corollary}

\numberwithin{equation}{section}
\newcommand{\bR}{\mathbf {R}}
\newcommand{\bv}{\mathbf {V}}
\newcommand{\bt}{\mathbf {\Theta}_n}
\newcommand{\ls}{\mathbf {\Lambda}_s}
\newcommand{\lb}{\mathbf {\Lambda}_b}
\newcommand{\bX}{\mathbf {X}}
\newcommand{\bx}{\mathbf {x}}
\newcommand{\bI}{\mathbf {I}}
\newcommand{\bus}{\mathbf {U}_s}
\newcommand{\bub}{\mathbf {U}_b}
\newcommand{\be}{\mathbf{e}}
\newcommand{\rlo}{\mathrm{log}}
\newcommand{\tr}{\operatorname{tr}}

\providecommand{\anon}{1}
\if1\anon
\hypersetup{pdfauthor={Tanyang Zhou; Yanqing Yin}}
\fi

\begin{document}

	\def\spacingset#1{\renewcommand{\baselinestretch}%
		{#1}\small\normalsize} \spacingset{1}

	\if1\anon
	{
		\renewcommand{\thefootnote}{\fnsymbol{footnote}}
		\begin{center}
{\LARGE\bfseries When Standardization Spreads Factor Variation Across the Spectrum\par}
\medskip
Tanyang Zhou$^{1}$ and Yanqing Yin$^{2}$\footnotemark[1]\\
$^{1}$School of Mathematics and Statistics, Chongqing University, Chongqing 401331, China\\
$^{2}$School of Statistics and Data Science, Nanjing Audit University, Nanjing 211815, China
\end{center}
\footnotetext[1]{Corresponding author: Yanqing Yin. Email: \href{mailto:yinyq799@nenu.edu.cn}{\nolinkurl{yinyq799@nenu.edu.cn}}.}
	} \fi
	
	\if0\anon
	{
		\bigskip
		\begin{center}
			{\LARGE\bf When Standardization Spreads Factor Variation Across the Spectrum}
		\end{center}
		\medskip
	} \fi
	
	\bigskip
	\begingroup
	\if1\anon\setlength{\leftmargini}{0pt}\fi
	\begin{abstract}
\spacingset{1.8}\normalsize
Standardizing variables by their sample standard deviations is routine in multivariate analysis and removes marginal scale differences before dependence is assessed. The resulting correlation matrix captures both dominant common components and dependence beyond them. A central inferential difficulty arises because the same random sample variances normalize the entire spectrum. In spiked covariance models, the sampling errors of leading and residual eigenvalues can become asymptotically independent. Sample standardization fundamentally changes this behavior and can preserve their dependence, allowing strong factors to remain a major source of uncertainty even after the leading eigenvalues are excluded. To obtain valid inference for both residual and full spectrum statistics, we establish their joint sampling law for sample correlation matrices with diverging spikes. The theory reveals exact trace compensation and statistic dependent transitions between asymptotic independence and dependence as the spikes grow. It leads to tests of an identity correlation matrix under strong factors and a feasible two sample comparison of residual spectral concentration. A sex balanced HapMap analysis illustrates the substantial effect of factor variation on residual spectrum inference.\end{abstract}
	
	\endgroup
	\noindent%
	\vfill
	
	\newpage
	\spacingset{1.8}
	
\section{Standardization and dependence inference}

Standardization is a routine first step in multivariate analysis when variables are measured on different scales. Dividing each variable by its sample standard deviation produces the sample correlation matrix, whose eigenvalues describe how standardized variation is distributed across orthogonal directions. In strongly correlated data, a few leading eigenvalues may reflect dominant common factors, while the residual spectrum characterizes dependence beyond these leading directions.

Previous work shows that sample standardization alters the second-order behavior of spiked eigenvalues and linear spectral statistics (LSS). Standardization precedes spectral decomposition, so sample variances inherit common-factor fluctuations and then rescale every coordinate. Spike and bulk errors may consequently reinforce or offset each other, allowing common-factor variation to remain relevant for residual-spectrum inference. A joint law can separate this source of sampling variation from changes in residual dependence.

For sample covariance matrices, spectral and LSS limits are well developed \citep{Silverstein1995,Bai2004,Pan2008,Pan2014,Zheng2015}, including LSS CLTs with an unbounded population spectral norm \citep{Yin2022Unbounded,Yin2026Liberating}. Limit laws are also available for the largest and spiked eigenvalues \citep{Johnstone2001,baik2006,paul2007,bai2008,bai2012,cai2020,jiang2021}. For their joint behavior, \citet{zhang2022} established asymptotic independence between suitably normalized spiked eigenvalues and LSS under a bounded population spectrum, while \citet{Liu2023} obtained asymptotically independent spike and bulk contributions for diverging spikes.

The corresponding theory for sample correlation matrices is less complete. LSS limits have been established under the identity population correlation matrix \citep{Gao2017} and for general population correlation matrices under Gaussian, independent-component, and elliptical models \citep{Mestre2017,yin2023}, with related results for rescaled correlation matrices and determinants \citep{Yin2021,Jiang2019}. For leading eigenvalues, \citet{Morales2021} identified normalization corrections in a block model, and \citet{r7} derived fluctuation limits under generalized spiked models with diverging spikes. The joint sampling behavior of the leading and residual spectra after standardization remains unresolved.

The distinction from covariance matrices appears in the reduction of the bulk statistic. In \citet{Liu2023}, the sum of the non-spiked eigenvalues reduces, up to a deterministic correction and a negligible remainder, to an LSS of a sample covariance matrix with bounded population spectrum. For correlation matrices, the analogous reduction retains the random normalization by sample variances of the full observations, including contributions from the strong factors. The resulting bulk statistic differs from the LSS of a correlation matrix constructed from residual data alone, and its marginal CLT requires analysis under the original normalization.

We derive a joint Gaussian law for the leading eigenvalues and residual spectral sums when a fixed number of population correlation spikes diverge, the bulk remains bounded, and the dimension and sample size grow proportionally. The theory covers real and proper complex observations with finite fourth moments, with the two components normalized at their respective nondegenerate scales. A diverging spectral gap can coexist with persistent spike--bulk dependence after standardization, whose form varies with the statistic and the spike growth rate.

Trace and logarithmic statistics illustrate two distinct mechanisms. Since a sample correlation matrix has unit diagonal, the centered sum of the leading eigenvalues is exactly the negative of the centered sum of the remaining eigenvalues at every spike strength. For the logarithmic statistic, Example~\ref{rem:decoupling} shows a transition from asymptotic independence to dependence as the spikes grow. The direct spike contribution to the centered full statistic vanishes across both regimes, while factor variation remains in the bulk variance once dependence emerges.

For residual dependence, consider comparing two genotype groups in which a few leading components capture major ancestry-related variation. After these components are removed, the remaining markers may still differ in how strongly their variation is concentrated across spectral directions. Two groups can share the same total factor strength and residual total variance while having different sums of squared residual eigenvalues. Our residual limit provides a feasible standard error for comparing this spectral concentration across groups.

A different question arises when the inferential target is overall dependence. Under the identity correlation null, strong-factor alternatives may concentrate evidence in one leading root or distribute comparable signal across several components. A largest-root test is naturally most sensitive to concentration in the leading direction, whereas logarithmic and quadratic spectral statistics accumulate evidence across multiple directions. The joint law gives their sampling distributions under strong-factor alternatives and permits direct comparison across different allocations of factor signal.

This paper is organized as follows. Section~\ref{Sec2} introduces the model and statistics, and Section~\ref{Sec3} establishes their joint sampling law. Section~\ref{app} develops the tests and residual-spectrum comparison, followed by simulation and real-data results in Sections~\ref{simula} and~\ref{sec:empirical}. Proofs and additional numerical results are collected in the Supplementary Material.

\section{Model and statistics for dependence analysis}\label{Sec2}

Random diagonal normalization is the source of the dependence described above. To isolate its role, we first express the sample correlation matrix in terms of population standardized observations and then introduce the statistics based on the full and residual spectra.

\subsection{Standardized observations and common components}

Let $\mathbf{y}:=\left(\mathbf{y}_1, \ldots, \mathbf{y}_n\right)$ be a random sample of size $n$ from a $p$-dimensional population with covariance matrix $\boldsymbol{\Sigma}_n$. Assume that $\mathbf{y}_j=\mathbf{\Gamma}_n \mathbf{x}_j, \quad j=1, \ldots, n,$ where $\boldsymbol{\Gamma}_n \boldsymbol{\Gamma}_n^*=\mathbf{\Sigma}_n$ and the entries of $\mathbf{x}_j=\left(x_{1 j}, \ldots, x_{p j}\right)^*$ are independently and identically distributed (i.i.d.), with mean zero and unit variance. All diagonal entries of $\boldsymbol{\Sigma}_n$ are strictly positive. Writing $\mathbf{\bX }_n=\left(\mathbf{x}_1, \ldots, \mathbf{x}_n\right)$, the sample covariance matrix is $\widehat{\mathbf{\mathbf{S}}}_n=\frac{1}{n} \sum_{j=1}^n \mathbf{y}_j \mathbf{y}_j^*=\frac{1}{n} \mathbf{Y}_n \mathbf{Y}_n^*$. The corresponding population and sample correlation matrices are, respectively,
$$
\mathbf{R}_n=\left[\operatorname{diag}\left(\boldsymbol{\Sigma}_n\right)\right]^{-1 / 2} \boldsymbol{\Sigma}_n\left[\operatorname{diag}\left(\boldsymbol{\Sigma}_n\right)\right]^{-1 / 2}, \quad \widehat{\mathbf{R}}_n=\left[\operatorname{diag}\left(\widehat{\mathbf{\mathbf{S}}}_n\right)\right]^{-1 / 2} \widehat{\mathbf{\mathbf{S}}}_n\left[\operatorname{diag}\left(\widehat{\mathbf{\mathbf{S}}}_n\right)\right]^{-1 / 2},
$$

Define $\mathbf{G}_n=\left(g_{k h}\right)_{k, h=1}^p=\left[\operatorname{diag}\left(\boldsymbol{\Sigma}_n\right)\right]^{-1 / 2} \boldsymbol{\Gamma}_n$, so that ${\mathbf{R}}_n=\mathbf{G}_n\mathbf{G}_n^*$. We suppress the subscript in $\bR$ and write the outer products of these standardized observations and their average as
$$\mathbf{\Theta}_j=\mathbf{G}_n \bx_j  \bx_j^* \mathbf{G}_n^*,  \quad \quad \mathbf{\Theta}_n=\frac{1}{n}\sum_{j=1}^n \mathbf{\Theta}_j =\frac{1}{n} \mathbf{G}_n \bX_n  \bX_n^* \mathbf{G}_n^*,$$
We have $\mathrm{E}(\mathbf{\Theta}_j)=\mathrm{E}(\mathbf{\Theta}_n)=\bR_n$ and
\begin{align*}
	\operatorname{diag}(\boldsymbol{\Theta}_n)
	&= \left[\operatorname{diag}(\boldsymbol{\Sigma}_n)\right]^{-1/2} \operatorname{diag}(\widehat{\mathbf{S}}_n) \left[\operatorname{diag}(\boldsymbol{\Sigma}_n)\right]^{-1/2}= \left[\operatorname{diag}(\boldsymbol{\Sigma}_n)\right]^{-1} \operatorname{diag}(\widehat{\mathbf{S}}_n),
\end{align*}

Thus, $	\widehat{\bR}_n$ can be written as $\widehat{\bR}_n = \left[\operatorname{diag}(\boldsymbol{\Theta}_n)\right]^{-1/2} \boldsymbol{\Theta}_n \left[\operatorname{diag}(\boldsymbol{\Theta}_n)\right]^{-1/2}.$ Let $\alpha_1,\alpha_2,\cdots,\alpha_p$ and $\lambda_1,\lambda_2,\cdots,\lambda_p$ denote the ordered eigenvalues of $\mathbf{R}_n$ and $\widehat{\bR}_n$, respectively. We partition the singular value decomposition of $\mathbf{G}_n$ into its spike and bulk components:
\begin{equation}
	\mathbf{G}_n=\mathbf{V \Lambda}^{1 / 2} \mathbf{U}^*=\left(\mathbf{V}_s, \mathbf{V}_b\right)\left(\begin{array}{cc}
		\mathbf{\Lambda}_s^{\frac{1}{2}} & 0 \\
		0 & \mathbf{\Lambda}_b^{\frac{1}{2}}
	\end{array}\right)\left(\mathbf{U}_s, \mathbf{U}_b\right)^* .
\end{equation}
Here $\mathbf{V}=(\bv_s,\bv_b)$ and $\mathbf{U}=(\mathbf{U}_s,\mathbf{U}_b)$ contain the left and right singular vectors of $\mathbf{G}_n$, respectively. Their first $M$ columns form the spike blocks, and the remaining $p-M$ columns form the bulk blocks, with $M$ fixed. The matrix $\mathbf{\Lambda}_s=\mathrm{diag}(\alpha_1,\alpha_2,...,\alpha_M)$ collects the simple spiked eigenvalues, which diverge, while $\mathbf{\Lambda}_b$ contains the bulk eigenvalues. For the bulk component, $\mathbf{G}_b=\bv_b \lb^{1 / 2} \mathbf{U}_b^*,\ \bR_b=\mathbf{G}_b\mathbf{G}_b^*$. Below, we suppress the subscript in $\mathbf{G}_n$.

We write $X_n\asymp Y_n$ when $X_n/C<Y_n<CX_n$ for some $C>1$, $X_n\backsimeq Y_n$ when $X_n=Y_n+o(Y_n)$, and $X_n=O(Y_n)$ when $|X_n|\leq C|Y_n|$. Here, the constant $C$ may vary from one occurrence to the next. For a matrix $\mathbf A_n$, $\|\mathbf A_n\|$ and $\mathrm{det}\mathbf A_n$ denote its spectral norm and determinant. We use $\be_k$ for the $k$th coordinate vector, $\boldsymbol\alpha^*$ for the conjugate transpose of $\boldsymbol\alpha$, and $\mathbf1_{\{X\geq x\}}$ for an indicator.

\subsection{Spectral summaries for factor and residual structure}

We distinguish two inferential targets, the full spectrum and the residual spectrum, while retaining the leading components to track the contribution of the dominant factors. For the ordered sample correlation eigenvalues $\lambda_1,\ldots,\lambda_p$, write $L_n(f)=\sum_{j=1}^p f(\lambda_j)$. The residual summary is based on $\sum_{j=M+1}^p f(\lambda_j)$, and the leading contribution $\sum_{h=1}^M f(\lambda_h)$ isolates the part associated with the estimated common components.

For the quadratic function, $\tr(\widehat{\bR}_n^2)-p =\sum_{i\ne j}|(\widehat{\bR}_n)_{ij}|^2$, so $L_n(x^2)-p$ measures aggregate squared sample correlation. The sum $\sum_{j=M+1}^p\lambda_j^2$ measures spectral concentration in the remaining directions. Before specifying their deterministic centers, we impose the following moment and dimensional conditions.

\begin{assumption}
	
	\label{ass1}
	The variables $x_{ij}$ are i.i.d., with $\mathrm{E}|x_{ij}|^4<\infty$ and common moments \label{A} $\mathrm{E}x_{ij}=0,\quad \mathrm{E}|x_{ij}|^2=1, \quad \beta_x=\mathrm{E}|x_{ij}|^4-|\mathrm{E}x_{ij}^2|^2-2, \quad \alpha_x=|\mathrm{E}x_{ij}^2|^2 .$ In the real case, $\mathbf{G}$ and the entries $x_{ij}$ are real and $\alpha_x=1$. In the proper complex case, $\mathrm{E}x_{ij}^2=0$, $\alpha_x=0$, and $\mathbf{G}$ may be complex.
	
\end{assumption}

\begin{assumption}
	\label{ass2}
	As $n \rightarrow \infty$, the ratio $c_n=p/n \rightarrow c \in (0,\infty) $.
\end{assumption}

\begin{assumption}
	\label{ass3}
	The bulk matrix $\mathbf{R_b=G_b G_b^*}$ has bounded spectral norm, and its empirical spectral distribution $H_n :=F^{\bR_b}$ converges weakly to a probability distribution $H$ on the real line.
\end{assumption}

\subsection{Deterministic centering and sampling scales}

Deterministic centering is described through empirical spectral distributions and their Stieltjes transforms. For a $p \times p$ Hermitian matrix $\mathbf{A}_n$ with eigenvalues ${\lambda_1 \geq \cdots \geq \lambda_p}$, its empirical spectral distribution (ESD) is $F^{\mathbf{A}_n}(x) = \frac{1}{p}\sum_{i=1}^p \mathbf{1} \left( \lambda_i(\mathbf{A}_n) \le x \right)$, where $\mathbf{1}(\cdot)$ is the indicator function. Throughout, we adopt the following sign convention for the Stieltjes transform.

\begin{definition}
	For any function $F$ of bounded variation on the real line, its Stieltjes transform is defined by $m_F(z) = \int \frac{1}{y - z} dF(y)$ for $z \in \mathbb{C}^+ \equiv \{z \in \mathbb{C}: \Im z > 0\}$.
\end{definition}

Under this model, the limiting ESD is characterized by the Mar\v{c}enko--Pastur equation of \citet{Silverstein1995}.

\begin{proposition}
	Under Assumptions \ref{ass1}--\ref{ass3}, $F^{\widehat{\mathbf{R}}_n}$ converges weakly almost surely to a nonrandom probability distribution $F^{c,H}$ as $n\to\infty$. For $z\in\mathbb C^+$, its Stieltjes transform ${m}_{F^{c,H}}={m}_{F^{c,H}}(z)$ satisfies
	\begin{equation}\label{mp}
		{m}_{F^{c,H}} = \int \frac{1}{t(1 - c - c z {m}_{F^{c,H}}) - z} dH(t),
	\end{equation}
	which is unique in the set $\left\{ {m}_{F^{c,H}} \in \mathbb{C} : -(1 - c)/z + c {m}_{F^{c,H}} \in \mathbb{C}^+ \right\}$. 
	
	Define the companion Stieltjes transform by $\underline{m}(z):=-\frac{1-c}{z}+c\,m_{F^{c,H}}(z)$. In terms of this transform, equation (\ref{mp}) becomes
	\begin{equation}\label{mp1}
		\underline{m} = -\left( z - c \int \frac{t}{1 + t \underline{m}} dH(t) \right)^{-1},
	\end{equation}
\end{proposition}

Solving for the spectral parameter gives the inverse relation
\begin{equation}\label{mp2}
	z = -\frac{1}{\underline{m}_{F^{c,H }}(z)}
	+ c \int \frac{t}{1 + t \underline{m}_{F^{c,H }}(z)} dH (t).
\end{equation}

For finite-$n$ residual centering, put $c_{nM}=(p-M)/n$ and $H_{2n}=F^{\lb}$. The bulk law $F^{c_{nM},H_{2n}}$ uses the retained $p-M$ coordinates, whereas $F^{c_n,H_n}$ uses their $p$-dimensional extension by $M$ zero eigenvalues. Since $H_n=(1-M/p)H_{2n}+(M/p)\delta_0$, the companion equations coincide and give the exact measure identity
$
pF^{c_n,H_n}=(p-M)F^{c_{nM},H_{2n}}+M\delta_0.
$
This identity accounts for the $Mf_\ell(0)$ correction in the full LSS below.

The bulk law provides the residual center, while the spike map determines the locations of the leading eigenvalues through 
\[
\varphi_{c,H}(z) = z\left(1 + \int \frac{ct}{z-t} dH(t)\right),
\quad
\psi_{c,H }(z) = \frac{\varphi_{c,H}(z)}{z} = 1 + \int \frac{c t}{z-t} dH(t).
\]
We write $\varphi=\varphi_{c,H}$ and $\varphi_n=\varphi_{c_n,H_n}$.
The derivative of the spike map,
$\varphi'_{c,H}(z)=1-c\int\frac{t^2}{(z-t)^2}\,dH(t)$,
appears in Assumption~\ref{ass4}. At an outlier location, the inverse relation further yields $1+\alpha_h \underline{m}_{F^{c,H}}(\varphi(\alpha_h))=0$.

\begin{assumption}
	\label{ass4}
	The spikes $\alpha_1,\dots,\alpha_M$ lie outside the support of $H$ and satisfy $\varphi'_{c,H}(\alpha_h) > 0,\  1\le h \le M$, and, when $M\geq2$, $\min_{j\neq k} \left| \frac{\alpha_k}{\alpha_j} - 1 \right| > d$, with $\min_{1\le h\le M}\alpha_h\to\infty$.
\end{assumption}

Residual fluctuations share a common scale determined by $\tr(\mathbf R^2)$, so the population spikes remain part of this scale even after the leading sample eigenvalues are excluded. Full spectral statistics carry an additional contribution from the derivative weights at the outlier locations. Write $\varphi_{n,h}=\varphi_n(\alpha_h)$ and $\boldsymbol\omega_{n,f}=(\omega_{n,f,h})_{h=1}^M{}^\top$, where $\omega_{n,f,h}=\varphi_{n,h}f'(\varphi_{n,h})/\sqrt n$, and define
\begin{equation}\label{An}
	s_{n,b}=\sqrt{\frac{n}{\tr(\mathbf R^2)}},\qquad
	s_{n,f}=\left[ \frac{\operatorname{tr}(\bR^2)}{n} + \frac{1}{n}\sum_{h=1}^{M} \left| \varphi_{n,h} f'(\varphi_{n,h}) \right|^2 \right]^{-1/2}=\dfrac{1}{\sqrt{s_{n,b}^{-2} + \|\boldsymbol \omega_{n,f}\|^2}}.
\end{equation}
Here $s_{n,b}$ provides the common residual scale, and $s_{n,f}$ combines it with the outlier derivative weights for the full statistic. By construction, $\left( \frac{s_{n,f}}{s_{n,b}} \right)^2 + \sum_{h=1}^M \big| s_{n,f} \omega_{n,f,h} \big|^2 = 1$. Consequently, $0 < s_{n,f} \leq s_{n,b}$ and $\max_{h\leq M} \big| s_{n,f} \omega_{n,f,h} \big| \leq 1$.

For a collection of test functions, write $s_{n,\ell}=s_{n,f_\ell}$. The corresponding renormalized LSSs are $\overline Y_\ell=s_{n,\ell}\mathcal Y_\ell$, where
\begin{align*}
	\mathcal Y_\ell
	&=\sum_{j=1}^p f_\ell(\lambda_j)-p\int f_\ell(x)\,dF^{c_n,H_n}(x)
	-\sum_{h=1}^M f_\ell(\varphi_n(\alpha_h))+Mf_\ell(0),
	\qquad \ell=1,\ldots,r.
\end{align*}

\subsection{Finite-$n$ decomposition of standardization effects}\label{sec:compact-kernels}

The spike and bulk blocks are subject to the same random normalization by sample variances, creating a direct channel through which their fluctuations can interact. A block decomposition makes this shared normalization explicit and leads to the resolvent and contrast kernels governing the LSS mean and covariance.

\subsubsection{A common diagonal normalization}

The block form below makes the shared diagonal rescaling explicit.
\[
\widehat{\bR}_n
=[\operatorname{diag}(\bt)]^{-1/2}\mathbf{V}
\underbrace{\begin{pmatrix}
		\mathbf{S}_{11}&\mathbf{S}_{12}\\
		\mathbf{S}_{21}&\mathbf{S}_{22}
\end{pmatrix}}_{\mathbf{S}}
\mathbf{V}^*[\operatorname{diag}(\bt)]^{-1/2},
\]
where the four blocks are
$
	\mathbf{S}_{11}
=\frac1n\ls^{1/2}\bus^*\bX\bX^*\bus\ls^{1/2},$ $
	\mathbf{S}_{12}=\mathbf{S}_{21}^*
=\frac1n\ls^{1/2}\bus^*\bX\bX^*\bub\lb^{1/2},$ $
	\mathbf{S}_{22}
=\frac1n\lb^{1/2}\bub^*\bX\bX^*\bub\lb^{1/2}.
$
The diagonal blocks $\mathbf{S}_{11}$ and $\mathbf{S}_{22}$ represent the spike and bulk directions, while $\mathbf{S}_{12}$ and $\mathbf{S}_{21}$ capture their interaction. All four blocks are rescaled by the same $\operatorname{diag}(\bt)$. This shared random normalization therefore provides a direct channel through which spike fluctuations enter statistics formed from the bulk eigenvalues.

\subsubsection{Bulk coordinates and contrasts for sample variances}

For the bulk component, $c_{nM}$ and $H_{2n}$ introduced above specify the finite-sample aspect ratio and population spectral distribution. We use $\circ$ for the Hadamard product, write $\mathbf A^{\circ2}=\mathbf A\circ\mathbf A$, and set $\mathbcal{L}_0(z)=\varphi'(z)/\psi(z)=z\varphi'(z)/\varphi(z)$. Variation in the sample variances is summarized by
\[
\bm{\mathfrak{R}}_{G}=(\alpha_x+1)(\bR\circ \overline{\bR})+\beta_x \left[(\mathbf{G}\circ \overline{\mathbf{G}})(\mathbf{G}^*\circ \mathbf{G}^\top)\right],
\]
where $m_{2n0}(z)$ is the Stieltjes transform of $F^{c_{nM},H_{2n}}$. Its companion transform is $m(z)=\underline m_{2n0}(z)=-\frac{1-c_{nM}}{z}+c_{nM}m_{2n0}(z)$. The matrices $\mathbf{K}_b=\mathbf{G}_b^*\mathbf{G}_b$ and $\boldsymbol{\Pi}_b=\mathbf{U}_b\mathbf{U}_b^* =\mathbf{I}-\mathbf{U}_s\mathbf{U}_s^*$ describe the retained bulk directions. We define
\[
\mathbf{J}_j(z)=\mathbf{U}_b\boldsymbol{\Lambda}_b
(\bI+\underline m_{2n0}(z)\boldsymbol\Lambda_b)^{-j}\mathbf{U}_b^*
=\mathbf{K}_b(\mathbf{I}+m(z)\mathbf{K}_b)^{-j},\qquad j=1,2.
\]
We also set $\Delta_n(z)=1-m(z)^2n^{-1}\tr(\mathbf{J}_1(z)^2)$. To express the effect of diagonal normalization in the original coordinates, define
\[
\mathbf{P}_j(z)=\mathbf{V}_b
(\bI+\underline m_{2n0}(z)\boldsymbol\Lambda_b)^{-j}\mathbf{V}_b^*
=\mathbf{V}_b\mathbf{V}_b^*(\mathbf{I}+m(z)\mathbf{R}_b)^{-j},\qquad j=1,2,
\qquad \mathbf{W}=\mathbf{G}\circ\overline{\mathbf{G}}.
\]
The identities
\[
\mathbf{G}\mathbf{J}_j(z)\mathbf{G}^*=\mathbf{R}_b^2\mathbf{P}_j(z),\qquad
\mathbf{P}_j(z)\mathbf{G}_b=\mathbf{G}_b(\mathbf{I}+m(z)\mathbf{K}_b)^{-j},\qquad
\mathbf{G}\mathbf{G}_b^*\mathbf{P}_j(z)=\mathbf{R}_b\mathbf{P}_j(z)
\]
connect these representations.
For each spike $h$, let $\mathbf{u}_h$ and $\mathbf{v}_h$ denote its right and left singular vectors, respectively, and set $\ell_h=\mathcal L_0(\alpha_h)$. The corresponding contrast matrices are
\[
\mathbf{H}_h=\mathbf{u}_h\mathbf{u}_h^*-\mathbf{G}^*[\mathbf{I}\circ(\mathbf{v}_h\mathbf{v}_h^*)]\mathbf{G},\qquad
\mathbf{H}_b(z)=\mathbf{J}_2(z)-\mathbf{G}^*[\mathbf{I}\circ\{\mathbf{R}_b\mathbf{P}_2(z)\}]\mathbf{G}.
\]
Each contrast separates the covariance contribution from the correction induced by normalization with the sample variances. The same matrices enter the spike, bulk, and mixed covariance terms below.

\subsubsection{Mean and covariance kernels}

The contrast matrices above determine four quantities in the joint law, namely the mean kernel, the spike covariance, the bulk covariance, and the covariance between spike and bulk fluctuations. All matrices below are evaluated at the displayed spectral parameter.

The unscaled mean kernel is
\begin{align*}
	B_n^0(z)&=\frac{m(z)^3}{n\Delta_n(z)}
	\left\{\frac{\alpha_x}{\Delta_n(z)}\tr(\mathbf{J}_1\mathbf{J}_2)
	+\beta_x\tr(\mathbf{J}_1\circ \mathbf{J}_2)\right\},\\
	b_n(z)&=B_n^0(z)-M\frac{m(z)+zm'(z)}{zm(z)}\\
	&\quad+\frac{m'(z)}{n}\Big[
	(1+\alpha_x)\tr\{(\mathbf{R}_b\mathbf{P}_1)\circ[\mathbf{R}_b(2\mathbf{P}_2-\mathbf{P}_1)]\}
	-\tr\{(\bm{\mathfrak R}_G\circ \mathbf{P}_1)\mathbf{R}_b\mathbf{P}_2\}\\
	&\qquad+\frac{\beta_x}{2}
	\tr\big\{[\mathbf{W}\circ((2\mathbf{P}_2-\mathbf{P}_1)\mathbf{G}_b)]\mathbf{G}_b^*\mathbf{P}_1
	+[\mathbf{W}\circ(\mathbf{P}_1\mathbf{G}_b)]\mathbf{G}_b^*(2\mathbf{P}_2-\mathbf{P}_1)\big\}\Big].
\end{align*}
Here $m'(z)=m(z)^2/\Delta_n(z)$.

Spike fluctuations are described by the covariance matrix in the following lemma.

\begin{lemma}[Simple-spike fluctuations]\label{lemma1}
	Under Assumptions~\ref{ass1}--\ref{ass4}, suppose that the limits defining $c_{ij}$ below exist and are finite. For any fixed subset of $K\leq M$ simple spikes, relabeled as $\alpha_1,\ldots,\alpha_K$, the vector $(\gamma_h)_{h=1}^K$, with $\gamma_h=\sqrt n\{\lambda_h(\widehat{\mathbf R}_n)/\varphi_n(\alpha_h)-1\}$, converges weakly to a mean-zero Gaussian vector with covariance matrix $\mathbf C=(c_{ij})$. In terms of the contrasts $\mathbf H_h$ defined above,
	\begin{align*}
		c_{ij}
		&=\lim_{n\to\infty}\ell_i\ell_j\Big[
		(1+\alpha_x)\tr(\mathbf H_i\mathbf H_j)
		+\beta_x\tr(\mathbf H_i\circ\mathbf H_j)\\
		&\hspace{9em}
		+(1+\alpha_x)\mathbf1_{\{i=j\}}
		\{\varphi'(\alpha_i)^{-1}-1\}\Big],
		\qquad 1\le i,j\le K.
	\end{align*}
\end{lemma}

For real observations, this lemma recovers the simple-spike CLT of \citet{r7}. The Supplementary Material provides the coordinate formulas and the extension to proper complex observations.

For the residual spectrum, the unscaled covariance kernel is
\begin{multline*}
	V_n(z_1,z_2)
	=(1+\alpha_x)
	\left\{\frac{m'(z_1)m'(z_2)}{[m(z_1)-m(z_2)]^2}
	-\frac{1}{(z_1-z_2)^2}\right\}\\
	+\frac{m'(z_1)m'(z_2)}{n}\Big[
	(1+\alpha_x)\{\tr(\mathbf{H}_b(z_1)\mathbf{H}_b(z_2))
	-\tr(\mathbf{J}_2(z_1)\mathbf{J}_2(z_2))\}+\beta_x\tr(\mathbf{H}_b(z_1)\circ \mathbf{H}_b(z_2))\Big].
\end{multline*}
Dependence between the spike and residual fluctuations is captured by the mixed kernel
\begin{equation*}
	\kappa_{n,h}^{(0)}(z)
	=-\frac{\ell_hm'(z)}{\sqrt n}
	\Big\{(1+\alpha_x)\tr(\mathbf{H}_h\mathbf{H}_b(z))
	+\beta_x\tr(\mathbf{H}_h\circ\mathbf{H}_b(z))-\beta_x\tr[(\mathbf{u}_h\mathbf{u}_h^*)\circ \mathbf{J}_2(z)]\Big\}.
\end{equation*}

\subsubsection{Contour forms of the LSS mean and covariance}

Contour integration against admissible analytic test functions converts the resolvent kernels into the mean and covariance terms of the centered LSSs. Choose contours in the common analytic domain that enclose the bulk support and exclude the sample outlier locations, with $\mathbcal C_1$ and $\mathbcal C_2$ disjoint. With the spike map, bulk law and mean kernel specified, the deterministically centered spike and residual contributions are
\[
\begin{aligned}
	S_n(f)&=\sum_{h=1}^M\{f(\lambda_h)-f(\varphi_n(\alpha_h))\},\\
	B_n(f)&=\sum_{j=M+1}^p f(\lambda_j)
	-(p-M)\int f(x)\,dF^{c_{nM},H_{2n}}(x)
	+\frac{1}{2\pi i}\oint_{\mathcal C}f(z)b_n(z)\,dz.
\end{aligned}
\]
Their sum is
\[
\begin{aligned}
	L_n(f)-\sum_{h=1}^M f(\varphi_n(\alpha_h))
	-(p-M)\int f(x)\,dF^{c_{nM},H_{2n}}(x)
	+\frac{1}{2\pi i}\oint_{\mathcal C}f(z)b_n(z)\,dz.
\end{aligned}
\]
Then, the full LSS can be written as
\[
\mathcal Y_\ell=S_n(f_\ell)+B_n(f_\ell)-\frac{1}{2\pi i}\oint_{\mathcal C}f_\ell(z)b_n(z)\,dz.
\]

For each test function, write $\mathcal Y_{n,\ell}=\mathcal Y_\ell$ and $\overline Y_{n,\ell}=\overline Y_\ell$, and define
\begin{gather*}
	\boldsymbol{\omega}_{n,\ell}=\boldsymbol{\omega}_{n,f_\ell}
	=\left(\frac{\varphi_n(\alpha_h)}{\sqrt n}
	f_\ell'(\varphi_n(\alpha_h))\right)_{h=1}^M{}^{\!\top},
	\quad \mu_{n,\ell}=-\frac{s_{n,\ell}}{2\pi i}
	\oint_{\mathbcal C}f_\ell(z)b_n(z)\,dz,\\
	\mathbf{k}_{n,\ell}=-\frac{1}{2\pi i}\oint_{\mathbcal C}
	f_\ell(z)(\kappa_{n,h}^{(0)}(z))_{h=1}^M\,dz,
	\quad B_{n,st}=-\frac{1}{4\pi^2}\oint_{\mathbcal C_1}\oint_{\mathbcal C_2}
	f_s(z_1)f_t(z_2)V_n(z_1,z_2)\,dz_1dz_2.
\end{gather*}
These quantities combine into the covariance matrix
\[
\Sigma_{n,st}=s_{n,s}s_{n,t}\left\{
\boldsymbol{\omega}_{n,s}^\top \mathbf{C}\boldsymbol{\omega}_{n,t}
+\boldsymbol{\omega}_{n,s}^\top \mathbf{k}_{n,t}
+\mathbf{k}_{n,s}^\top\boldsymbol{\omega}_{n,t}+B_{n,st}\right\}.
\]
\section{Joint sampling distributions}\label{Sec3}

The kernels derived above determine the joint fluctuations of the leading eigenvalues and residual spectral sums, and hence the limiting laws of statistics based on the full spectrum.

\subsection{Joint limit theorems}
To transfer the outlier fluctuations to an LSS, we linearize each test function around the corresponding deterministic outlier location $\varphi_n(\alpha_h)$. The following condition ensures that the derivative remains asymptotically unchanged over the relative fluctuation scale of the sample outlier.

\begin{assumption}
	The test functions $f_1,\dots,f_r$ are analytic on a connected open region containing the support of $F^{c_n,H_n}$ and the deterministic outlier locations $\varphi_n(\alpha_h)$, $1\le h\le M$, for almost all $n$, and satisfy $f_\ell(x)\in\mathbb R$ for real $x$ in this region. For each nonconstant $f_\ell$, its derivative is defined and nonzero for all sufficiently large positive real arguments, and
	\[
	\lim_{\substack{x_n,y_n\to+\infty\\x_n/y_n\to1}}
	\frac{f_\ell'(x_n)}{f_\ell'(y_n)}=1.
	\]
	\label{ass5}
\end{assumption}
At the component level, we collect the relative spike errors $\gamma_h$ from Lemma~\ref{lemma1} and the residual spectral sums normalized by the common scale $s_{n,b}$ in \eqref{An}.
\[
\mathbf S_n^{\mathrm{sp}}=(\gamma_1,\ldots,\gamma_M)^\top,\qquad
\mathbf S_n^{\mathrm{bu}}=s_{n,b}\bigl(B_n(f_1),\ldots,B_n(f_r)\bigr)^\top.
\]
The kernels in Section~\ref{sec:compact-kernels} determine the corresponding covariance blocks
\[
	\boldsymbol\Omega_n^{\mathrm{ss}}=\mathbf C,\quad
	\boldsymbol\Omega_n^{\mathrm{sb}}=s_{n,b}(\mathbf k_{n,1},\ldots,\mathbf k_{n,r}),\quad
	\boldsymbol\Omega_n^{\mathrm{bs}}=(\boldsymbol\Omega_n^{\mathrm{sb}})^\top,\quad
	\boldsymbol\Omega_n^{\mathrm{bb}}=\bigl(s_{n,b}^2B_{n,st}\bigr)_{s,t=1}^r,
\]
and set
\[
\boldsymbol\Omega_n=
\begin{pmatrix}
	\boldsymbol\Omega_n^{\mathrm{ss}}&\boldsymbol\Omega_n^{\mathrm{sb}}\\
	\boldsymbol\Omega_n^{\mathrm{bs}}&\boldsymbol\Omega_n^{\mathrm{bb}}
\end{pmatrix}.
\]

\begin{theorem}[Joint limit for spike and residual fluctuations]\label{th:joint-components}
	Under Assumptions~\ref{ass1}--\ref{ass5}, let $r$ and $M$ be fixed. Suppose that $\boldsymbol\Omega_n\to\boldsymbol\Omega$ for a finite matrix $\boldsymbol\Omega$. Then
	\[
	\begin{pmatrix}
		\mathbf S_n^{\mathrm{sp}}\\
		\mathbf S_n^{\mathrm{bu}}
	\end{pmatrix}
	\xrightarrow{d}
	N_{M+r}(\mathbf 0,\boldsymbol\Omega).
	\]
\end{theorem}
 In the contrast $\mathbf H_b(z)=\mathbf J_2(z)-\mathbf G^*[\mathbf I\circ\{\mathbf R_b\mathbf P_2(z)\}]\mathbf G$, the matrices $\mathbf J_2(z)$ and $\mathbf R_b\mathbf P_2(z)$ describe the bulk structure, while the outer factors $\mathbf G$ retain the full population structure. Diverging spikes can therefore affect the residual covariance even after their sample eigenvalues are removed. Their influence also enters through the scale $s_{n,b}=\{n/\tr(\mathbf R^2)\}^{1/2}$, which depends on the full population correlation matrix. Dependence between the two components is encoded by the off-diagonal blocks $\boldsymbol\Omega^{\mathrm{sb}}$ and $\boldsymbol\Omega^{\mathrm{bs}}$. A nonzero limiting block gives asymptotic dependence between the spike and residual fluctuations, so a diverging spectral gap can coexist with dependence at their respective sampling scales.

Statistics based on the full spectrum combine these two sources of fluctuation with different weights. Under Assumption~\ref{ass5},
\[
s_{n,\ell}\mathcal Y_{n,\ell}-\mu_{n,\ell}
=s_{n,\ell}\boldsymbol\omega_{n,\ell}^\top\mathbf S_n^{\mathrm{sp}}
+\frac{s_{n,\ell}}{s_{n,b}}(\mathbf S_n^{\mathrm{bu}})_\ell+o_p(1).
\]
This decomposition transfers the component level joint law to the full spectral statistics. We first state the scalar limit.

\begin{theorem}
	\label{th1}
	Under Assumptions \ref{ass1}--\ref{ass5}, suppose
	$\liminf_{n\to\infty}\Sigma_{n,11}>0$. With the quantities defined in Section~\ref{sec:compact-kernels},
	\[
	\frac{\overline Y_1-\mu_{n1}}{\sigma_{n1}}
	\xrightarrow{d}N(0,1),\qquad \sigma_{n1}^2=\Sigma_{n,11}.
	\]
	In particular,
	\begin{equation}\label{sigma}
		\sigma_{n1}^2=s_{n,1}^2
		\left\{\boldsymbol{\omega}_{n,1}^\top \mathbf{C}\boldsymbol{\omega}_{n,1}
		+2\boldsymbol{\omega}_{n,1}^\top \mathbf{k}_{n,1}+B_{n,11}\right\}
		=\sigma_{\mathrm{spi}}^2+\sigma_{\mathrm{spi,bulk}}^2+\sigma_{\mathrm{bulk}}^2.
	\end{equation}
\end{theorem}

The multivariate form follows by collecting the statistics for $f_1,\ldots,f_r$.

\begin{theorem}
	\label{th2}
	Under Assumptions~\ref{ass1}--\ref{ass5}, let $r$ and $M$ be fixed. Suppose that $\boldsymbol\Sigma_n\to\boldsymbol\Sigma$ for a finite matrix $\boldsymbol\Sigma$. Let
	$\overline Y_{n,\ell}=s_{n,\ell}\mathcal Y_{n,\ell}$ for $1\le\ell\le r$, and let
	$\boldsymbol{\mu}_n=(\mu_{n,1},\ldots,\mu_{n,r})^\top$ and $\boldsymbol{\Sigma}_n=(\Sigma_{n,st})$
	be defined in Section~\ref{sec:compact-kernels}. Then
	\[
	\overline{\mathbf{Y}}_n-\boldsymbol{\mu}_n\xrightarrow{d}N_r(0,\boldsymbol{\Sigma}).
	\]
\end{theorem}

\subsection{Factor strength and sampling scales}
Equation~\eqref{sigma} makes the role of factor strength explicit. For $f(x)=x$, the identity $\tr(\widehat{\bR}_n)=p$ forces the random parts of $S_n(x)$ and $B_n(x)$ to cancel exactly. The resulting full statistic is deterministic and falls outside the nondegenerate setting of Theorem~\ref{th1}.

For logarithmic and power functions, the contribution from diverging spikes depends on the derivative weights at their sample locations. We compare $\max_{1\le h\le M}|\alpha_h f'(\varphi_n(\alpha_h))|$ with $\sqrt{\tr(\bR^2)}\asymp\sqrt{n+\sum_{h=1}^M\alpha_h^2}$ and define
\[
\theta_{n,f} = \frac{\sum_{h=1}^M \left|\varphi_n(\alpha_h)f'(\varphi_n(\alpha_h))\right|^2}{\mathrm{tr}(\bR^2)}.
\]
The ratio $\theta_{n,f}$ summarizes the magnitude of the outlier derivative weights relative to the bulk contribution. The variance itself is determined by the full expression in \eqref{sigma}, including possible cancellation or degeneracy among its terms.

\begin{example}[Logarithmic and power statistics]
	\leavevmode\par\nobreak
	\smallskip\noindent\textit{Logarithmic statistic.}
	For $f(x)=\rlo(x)$ and diverging spikes, $s_{n,f}$ and $s_{n,b}$ are asymptotically equivalent. Since $\alpha_h f'(\alpha_h)=1$, the direct spike contribution vanishes. Diagonal normalization nevertheless contributes to the bulk variance. When $\alpha_1=o(\sqrt n)$, the unscaled kernel $V_n(z_1,z_2)$ contains the contribution
	\begin{align*}
		\frac{m'(z_1)m'(z_2)}{n}\Big[
		&(1+\alpha_x)\tr\{[\mathbf{J}_2(z_1)-\mathbf{H}_b(z_1)]
		[\mathbf{J}_2(z_2)-\mathbf{H}_b(z_2)]\}\\
		&+\beta_x\tr\{[\mathbf{J}_2(z_1)-\mathbf{H}_b(z_1)]
		\circ[\mathbf{J}_2(z_2)-\mathbf{H}_b(z_2)]\}\Big].
	\end{align*}
	
	\smallskip\noindent\textit{Quadratic statistic.}
	When $f(x)=x^2$ and $\alpha_1\to\infty$, $s_{n,f} \asymp \left(1+\frac{\alpha_1^4}{n}\right)^{-1/2}$.
	The largest derivative weight is $2\alpha_1^2$, which reaches the bulk scale at $\alpha_1\asymp n^{1/4}$. The spike coefficient vanishes when $\alpha_1=o(n^{1/4})$ and becomes dominant in this scale comparison when $n^{1/4}=o(\alpha_1)$.
	
	\smallskip\noindent\textit{Higher powers.}
	For $f(x)=x^m$, $m>2$, the largest derivative weight is $m\alpha_1^m$. It reaches the bulk scale at $\alpha_1\asymp n^{1/(2m)}$. The spike coefficient vanishes below this threshold and dominates the corresponding scale comparison above it.
	
\end{example}

\subsection{Dependence between spike and residual fluctuations}
These thresholds describe the relative contributions to statistics based on the full spectrum. Dependence between the spike and residual fluctuations is governed separately by the mixed covariance kernel. A nonzero limiting mixed covariance produces asymptotic dependence between the two components. In the sample covariance setting, \citet{Liu2023} obtained asymptotically independent limits for the corresponding quantities.

The trace exhibits an exact form of dependence. Whenever the component variance is positive, the compensation above gives perfect negative correlation at every spike strength. More generally, for any constants $a,b$, $L_n(f+ax+b)=L_n(f)+(a+b)p$, while $S_n(f+ax)=S_n(f)+aS_n(x)$ and $B_n(f+ax)=B_n(f)+aB_n(x)$. Hence an affine change of the test function preserves the random fluctuation of the full spectral statistic while redistributing the trace fluctuation between its two components.

\begin{example}[Dependence at separate scales]\label{rem:decoupling}
	\leavevmode\par\nobreak
	\smallskip\noindent\textit{Gaussian logarithmic model.}
	Consider real Gaussian observations with
	\[
	\bR=(1-\rho_n)\bI_p+\rho_n\mathbf{1}_p\mathbf{1}_p^\top,\qquad
	\alpha_1=1+(p-1)\rho_n\longrightarrow\infty,\qquad
	\frac{\alpha_1}{\sqrt n}\longrightarrow\tau\in[0,\infty),
	\]
	where $p/n\to c\in(0,1)$ and $\mathbf G=\bR^{1/2}$. The observations admit the factor representation $y_{ij}=\sqrt{\rho_n}\,F_j+\sqrt{1-\rho_n}\,\varepsilon_{ij}, \qquad 1\le i\le p,\quad 1\le j\le n$, where $F_j$ and $\varepsilon_{ij}$ are independent standard normal variables. Each population marginal variance is one, while every sample variance contains the common centered factor term $\rho_n n^{-1}\sum_{j=1}^n(F_j^2-1)$.

	Put $\delta_e=1-\rho_n$, $q_n=1-1/p$ and
	$\mathbf P_0=\bI_p-p^{-1}\mathbf{1}_p\mathbf{1}_p^\top$.
	For this model, the spike correction matrix and its logarithmic bulk counterpart are
	\[
	\mathbf H_1=p^{-1}(\mathbf{1}_p\mathbf{1}_p^\top-\bR),\qquad
	\mathbf H_L=\mathbf P_0-q_n\bR=-\frac{\alpha_1}{\delta_e}\mathbf H_1.
	\]
	Hence $\tr(\mathbf H_1^2)=q_n\delta_e^2,\qquad \tr(\mathbf H_1\mathbf H_L)=-q_n\delta_e\alpha_1,\qquad \tr(\mathbf H_L^2)=q_n\alpha_1^2$. Substitution into the formulas for equal bulk eigenvalues in Section~\ref{sec:classical-tests} gives
	\[
	e_\bR(1)=-\frac{2\ell_1q_n\delta_e\alpha_1}{\sqrt n},\qquad
	B_{n,LL}=-2\{\log(1-c_{nM})+c_{nM}\}
	+\frac{2q_n\alpha_1^2}{n}.
	\]
	Since $\delta_e\to1$ and $\ell_1\to1$, the spike variance, mixed covariance and bulk variance converge to $2$, $-2\tau$ and $2\{\tau^2-\log(1-c)-c\}$, respectively. Theorem~\ref{th:joint-components} then yields
	\[
	\begin{pmatrix}
		\sqrt n\,S_n(\log)\\[2pt]
		B_n(\log)
	\end{pmatrix}
	\xrightarrow{d}
	N_2\left(
	\mathbf 0,\,
	2\begin{pmatrix}
		1&-\tau\\
		-\tau&\tau^2-\log(1-c)-c
	\end{pmatrix}\right).
	\]
	The limiting correlation is $-\tau/\sqrt{\tau^2-\log(1-c)-c}$, which vanishes at $\tau=0$ and is nonzero for $\tau>0$. Thus $\alpha_1\asymp\sqrt n$ marks the transition from asymptotic independence to dependence when the two components are viewed at their own fluctuation scales.

	The logarithmic spike term has variance of order $n^{-1}$, so its direct contribution and the mixed contribution to the full logarithmic LSS variance vanish throughout this regime. Let $Z$ and $W$ be independent centered Gaussian variables with variances $2$ and $2\{-\log(1-c)-c\}$, respectively. Equivalently, $\bigl(\sqrt n\,S_n(\log),\,B_n(\log)\bigr) \xrightarrow{d}(Z,\,W-\tau Z)$, and $S_n(x)=\tau\sqrt n\,S_n(\log)+o_p(1)$. Together with $B_n(x)=-S_n(x)+o_p(1)$, this gives, for $\tau>0$, $\bigl(S_n(\log-x),\,B_n(\log-x)\bigr) \xrightarrow{d}(-\tau Z,\,W)$. Subtracting $x$ therefore transfers the trace fluctuation between the two components while preserving the limiting centered statistic based on the full spectrum.
	
	\smallskip\noindent\textit{Quadratic statistic.}
	For $f(x)=x^2$, substituting $e_\bR(1)$ into the mixed term for equal bulk eigenvalues gives, throughout $\alpha_1=O(\sqrt n)$,
	\[
	\sigma_{\text{spi,bulk}}^2 = -16 s_{n,W}^2 \ell_1 q_n \delta_e^3 \big(1 + c_{nM}\big) \frac{\varphi_n(\alpha_1)^2 \alpha_1}{n}=o(1).
	\]
	The mixed contribution is therefore governed by the overlap $e_\bR(1)$, in particular $\tr(\mathbf H_1\mathbf H_L)$, together with the remaining factors in the covariance term. Spike strength alone does not determine its order.
\end{example}

The logarithmic and quadratic rates describe different features of the joint fluctuation. For $f(x)=\log x$, $\alpha_1\asymp\sqrt n$ marks the dependence threshold for the leading and residual components at their respective fluctuation scales. For $f(x)=x^2$, $\alpha_1\asymp n^{1/4}$ marks the scale at which the spike derivative weight matches the bulk contribution to the full statistic, while the mixed covariance above remains $o(1)$.

\section{Inference on residual and overall dependence}\label{app}
Joint sampling laws yield inference on both residual and overall dependence. The residual limit yields a feasible two-sample comparison of spectral concentration after the leading components are removed, with a standard error that retains the contribution of strong factors. For overall dependence, logarithmic and quadratic statistics based on the full spectrum lead to tests of an identity correlation matrix and allow the effect of different allocations of factor signal on power to be examined.
\subsection{Two-sample comparison of residual structure}\label{sec:reference-bulk}

At a common total factor strength, residual spectral concentration measures how dependence is distributed beyond the dominant components. Such residual structure is relevant to covariance estimation and portfolio allocation \citep{FanLiaoMincheva2013}, while related detection problems arise in latent factor models \citep{PesaranXie2026}. 

Let $\mathbf R_0$ and $\mathbf R_1$ be the population correlation matrices of two independent samples, each of size $n$. Observations in group $\ell$ have the representation $\mathbf y_{\ell r}=\mathbf G_\ell\mathbf x_{\ell r}$, $1\le r\le n$, with $\mathbf G_\ell\mathbf G_\ell^\top=\mathbf R_\ell$. The coordinate distributions in each group satisfy the real-valued case of Assumption~\ref{ass1}, with fourth cumulants $\beta_\ell=\mathrm{E}x_{\ell,11}^4-3$. Fix $M\geq1$, put $q=p-M$ and let $p/n\to c\in(0,\infty)$. The spikes satisfy $\alpha_{\ell,h}\asymp p$, $1\leq h\leq M$, $\ell=0,1$, while the residual eigenvalues obey
\begin{equation}\label{eq:reference-bulk-model}
	t_{\ell,j}=\delta_e+\varepsilon_{\ell,j,n},\qquad
	\sum_{j=1}^{q}\varepsilon_{\ell,j,n}=0,\qquad
	\varepsilon_n^\star:=\max_{\ell\in\{0,1\}}\max_{j\leq q}
	|\varepsilon_{\ell,j,n}|\longrightarrow0,\qquad \ell=0,1,
\end{equation}
where $\delta_e\in(0,1)$ is fixed. The unit diagonal constraint gives $\sum_{h=1}^{M}\alpha_{\ell,h}=p-q\delta_e$, so the two groups have the same total spike strength while allowing different allocations and eigenspaces.

At the common residual total variance, we measure the change in residual spectral concentration by
\begin{equation}\label{eq:reference-target}
	\Delta_{b,n}:=\sum_{j=1}^{q}t_{1,j}^2-\sum_{j=1}^{q}t_{0,j}^2
	=\sum_{j=1}^{q}\varepsilon_{1,j,n}^2
	-\sum_{j=1}^{q}\varepsilon_{0,j,n}^2.
\end{equation}
Because both groups have residual total variance $q\delta_e$, a positive $\Delta_{b,n}$ corresponds to a smaller residual effective rank in group $1$, where $r_{\mathrm{eff},\ell}=(q\delta_e)^2/\sum_jt_{\ell,j}^2$. Under \eqref{eq:reference-bulk-model}, we test $H_{0,b}:\Delta_{b,n}=0$ against $H_{1,b}:\Delta_{b,n}>0$. The null equates residual spectral concentration while allowing the covariance matrices and residual spectra to differ. Quadratic-form tests such as \citet{LiChen2012} aggregate dense matrix differences, while entrywise procedures such as \citet{CaiLiuXia2013} target sparse changes. Our target is the second residual spectral moment and is invariant to eigenvector changes and spectral changes that preserve this moment.

Residual quadratic sums provide a natural estimator of this contrast. Writing $\lambda_{\ell,1}\geq\cdots\geq\lambda_{\ell,p}$ for the ordered sample correlation eigenvalues, define
\begin{equation}\label{eq:reference-statistic}
	Q_{\ell,n}=\sum_{j=M+1}^{p}\lambda_{\ell,j}^2,\qquad
	D_{b,n}=Q_{1,n}-Q_{0,n},\qquad
	\widehat\delta_{\ell,n}=\frac1q\sum_{j=M+1}^{p}\lambda_{\ell,j}.
\end{equation}
Applying Theorem~\ref{th:joint-components} with $f(x)=x^2$ gives the bulk center for group $\ell$, $m_{\ell,n}:=\sum_{j=1}^q t_{\ell,j}^2 +\frac1n\left(\sum_{j=1}^q t_{\ell,j}\right)^2 =\sum_{j=1}^q t_{\ell,j}^2+\frac{(q\delta_e)^2}{n},$ since $\sum_jt_{\ell,j}=q\delta_e$. Equal sample sizes therefore imply $m_{1,n}-m_{0,n}=\Delta_{b,n}$.

In each group, the leading variance is governed by fluctuations in the squared length of an observation. Put $\mathbf K_\ell=\mathbf G_\ell^\top\mathbf G_\ell$ and
\begin{equation}\label{eq:reference-energy-variance}
	\nu_{\ell,n}:=\operatorname{Var}(\|\mathbf y_{\ell1}\|^2)
	=2\tr(\mathbf R_\ell^2)+\beta_\ell\tr(\mathbf K_\ell\circ\mathbf K_\ell).
\end{equation}
For Gaussian coordinates, $\nu_{\ell,n}=2\sum_{h\le M}\alpha_{\ell,h}^2+O(p)$. The fourth-cumulant term further captures the contribution of the right singular vectors of $\mathbf G_\ell$. With $c_{nM}=q/n$, the resulting variance under strong factors is
\begin{equation}\label{eq:reference-variance}
	V_{\ell,n}=\frac{4\delta_e^4(1+c_{nM})^2}{n}\nu_{\ell,n},
	\qquad V_n=V_{0,n}+V_{1,n}.
\end{equation}

For group $\ell$, the exact unscaled bulk variance $B_{\ell,n,WW}$ is obtained from the covariance kernel in Section~\ref{sec:compact-kernels} with $f(x)=x^2$ and reduces to \eqref{eq:flat-bulk-variance-W} when the bulk eigenvalues are equal. Under \eqref{eq:reference-bulk-model} and the nondegeneracy condition in Corollary~\ref{prop:reference-bulk}, Section~S.7 of the Supplementary Material establishes $B_{\ell,n,WW}=V_{\ell,n}+o(p)$ and $V_{\ell,n}\asymp p$ for $\ell=0,1$. This representation leads to a feasible variance estimator based on standardized squared lengths. For $\mathbf Y_\ell=(\mathbf y_{\ell1},\ldots,\mathbf y_{\ell n})$ and $\mathbf S_\ell=n^{-1}\mathbf Y_\ell\mathbf Y_\ell^\top$, define
\begin{align}\label{eq:reference-energy-estimator}
	e_{\ell r}&=\mathbf y_{\ell r}^\top
	[\operatorname{diag}(\mathbf S_\ell)]^{-1}\mathbf y_{\ell r},
	&\widehat\nu_{\ell,n}&=\frac1n\sum_{r=1}^n(e_{\ell r}-p)^2,\\
	\widehat V_{\ell,n}&=\frac{4\widehat\delta_{\ell,n}^4(1+c_{nM})^2}{n}
	\widehat\nu_{\ell,n},
	&\widehat V_n&=\widehat V_{0,n}+\widehat V_{1,n}.
	\label{eq:reference-variance-estimator}
\end{align}
Because the empirical mean of $e_{\ell r}$ is exactly $p$, $\widehat\nu_{\ell,n}$ is their empirical variance. The same diagonal standardization used to form the sample correlation matrix therefore also provides the variance estimate for the residual comparison.

\begin{corollary}[Two-sample comparison under strong factors]\label{prop:reference-bulk}
	Under Assumptions~\ref{ass1}--\ref{ass4} and the real-valued model and spectral conditions of this subsection, suppose $\liminf_{n\to\infty}\nu_{\ell,n}/p^2>0$ for $\ell=0,1$. Then
	\begin{equation}\label{eq:reference-clt}
		\frac{D_{b,n}-\Delta_{b,n}}{\sqrt{V_n}}\xrightarrow{d}N(0,1),
		\qquad \frac{\widehat V_n}{V_n}\xrightarrow{p}1.
	\end{equation}
	Consequently, $Z_{b,n}:=D_{b,n}/\sqrt{\widehat V_n}$ converges to $N(0,1)$ under $H_{0,b}$. If $\Delta_{b,n}/\sqrt{V_n}\to d\in[0,\infty)$, rejection at $Z_{b,n}>z_{1-\gamma}$ has limiting power
	\begin{equation}\label{eq:reference-power}
		1-\Phi(z_{1-\gamma}-d),
	\end{equation}
	where $\Phi$ is the standard normal distribution function and $z_{1-\gamma}=\Phi^{-1}(1-\gamma)$. The power converges to one when $\Delta_{b,n}/\sqrt{V_n}\to\infty$.
\end{corollary}

The equal-bulk case $\varepsilon_{\ell,j,n}=0$ makes the contribution of the strong factors to the standard error explicit. A benchmark based only on the residual correlation spectrum treats the eigenvalues divided by $\delta_e$ as those of a $q$-dimensional sample correlation matrix with identity population correlation. Under Gaussian coordinates, the quadratic LSS has variance $4c_{nM}^2$ in each group, giving
\begin{equation}\label{eq:reference-bulk-only-variance}
	V_n^{\mathrm{bulk}}=8\delta_e^4c_{nM}^2,\qquad
	\frac{V_n}{V_n^{\mathrm{bulk}}}
	=\frac{(1+c_{nM})^2}{nc_{nM}^2}
	\left(\sum_{\ell=0}^1\sum_{h=1}^M\alpha_{\ell,h}^2+2q\delta_e^2\right)\asymp p.
\end{equation}
Thus the residual-spectrum benchmark understates the sampling variance by an order of $p$. For any $\gamma<1/2$, the nominal level-$\gamma$ test based on $Z'_{b,n}=\frac{D_{b,n}}{\sqrt{V_n^{\mathrm{bulk}}}}$ satisfies
\[
\Pr_{H_{0,b}}\left(\frac{D_{b,n}}{\sqrt{V_n^{\mathrm{bulk}}}}>z_{1-\gamma}\right)\longrightarrow\frac12.
\]
Factor variation therefore remains essential for valid inference even after the leading sample eigenvalues are excluded.

Corollary~\ref{prop:reference-bulk} also yields the two-sided confidence interval $D_{b,n}\pm z_{1-\gamma/2}\sqrt{\widehat V_n}$ for $\Delta_{b,n}$. If the squared perturbation energies of the two groups differ by order $pa_n^2$, then $\Delta_{b,n}\asymp pa_n^2$ and $\sqrt{V_n}\asymp p^{1/2}$, giving the local scale $a_n\asymp p^{-1/4}$.

\subsection{Testing an identity correlation matrix}\label{sec:classical-tests}

We next use the full spectrum to test an identity correlation matrix under alternatives with a finite number of diverging spikes and equal bulk eigenvalues.
\begin{equation}\label{t1}
	H_0: \bR = \mathbf{I}_p \quad \text{vs.} \quad H_1: \bR = \mathbf{V}
	\begin{pmatrix}
		\ls & \mathbf{0} \\
		\mathbf{0} & \delta_e \mathbf{I}_{p-M}
	\end{pmatrix}
	\mathbf{V}^*,
\end{equation}
where $\ls$ contains the diverging spikes and $\delta_e=\delta_{e,n}\in(0,1)$ is the common bulk eigenvalue, with $\liminf_{n\to\infty}\delta_e>0$. The unit diagonal condition $\bR_{ii}=1$ for $i=1,\cdots,p$ implies $\sum_{h=1}^M\alpha_h+(p-M)\delta_e=p$. Define the bulk leverage coefficient $\vartheta_{n,b}=(p-M)^{-1}\tr(\boldsymbol{\Pi}_b^{\circ2})$, which enters the fourth-cumulant terms through $\beta_x\vartheta_{n,b}$. Under $H_0$, we take $\mathbf G_n=\mathbf I_p$.

Following \citet{yin2023}, consider the logarithmic statistic $T_L=\log|\widehat{\bR}_n|=\sum_{j=1}^p\log\lambda_j$ and Nagao's quadratic statistic $T'_W=\tr(\widehat{\bR}_n-\bI_p)^2$. The identity $\tr(\widehat{\bR}_n)=p$ gives $T'_W=T_W-p$, where $T_W=\tr(\widehat{\bR}_n^2)=\sum_{j=1}^p\lambda_j^2$, so the quadratic analysis can be based on $T_W$.

For $f_L(x)=\log x$, the contour is chosen within positive neighborhoods of the bulk spectrum and outlier locations. Since $\lambda_{\min}(\widehat{\bR}_n)\ge\lambda_{\min}(\bR)\lambda_{\min}(\bX\bX^*/n)/\|\operatorname{diag}(\boldsymbol\Theta_n)\|>0$ with probability tending to one, define $T_L=0$ on the negligible singular event. Theorem~\ref{th1} applies directly to $f_W(x)=x^2$. With $s_{n,a}=s_{n,f_a}$ for $a\in\{L,W\}$, the sampling scales are $s_{n,L} = \sqrt{\frac{n}{\operatorname{tr}(\bR^2)+M}},\quad s_{n,W} = \sqrt{\frac{n}{\operatorname{tr}(\bR^2)+4\sum_h \varphi_{n,h}^4}}.$

Equal bulk eigenvalues simplify the centering and covariance terms in Section~\ref{sec:compact-kernels}. Put
\begin{equation}\label{eq:flat-bulk-HL}
	\mathbf P_0=\mathbf V_b\mathbf V_b^*,\qquad
	\mathbf H_L=\boldsymbol{\Pi}_b-
	\mathbf G^*(\mathbf I\circ\mathbf P_0)\mathbf G.
\end{equation}
The logarithmic and quadratic statistics have spike weights
\begin{equation}\label{eq:flat-bulk-weights}
	\boldsymbol{\omega}_{n,L}=\frac{\mathbf 1_M}{\sqrt n},\qquad
	\boldsymbol{\omega}_{n,W}
	=\frac{2}{\sqrt n}\bigl(\varphi_n^2(\alpha_h)\bigr)_{h=1}^M.
\end{equation}
Their mean corrections involve
\begin{align}\label{eq:flat-bulk-mean-terms}
	b_\bR
	&=\frac{1}{2n}\tr\{(\bm{\mathfrak R}_G\circ\mathbf P_0)\mathbf P_0\},\\
	a_\bR
	&=\frac{(1+\alpha_x)\delta_e}{n}\tr(\mathbf P_0^{\circ2})
	+\frac{\beta_x}{n}
	\tr\{\mathbf W[(\mathbf U_b\mathbf V_b^*)
	\circ\overline{\mathbf U_b\mathbf V_b^*}]\}.
\end{align}
Dependence between the leading and bulk fluctuations is represented by
\begin{align}\label{eq:flat-bulk-mixed}
	e_\bR(h)
	&=\frac{\ell_h}{\sqrt n}\Big[
	(1+\alpha_x)\tr(\mathbf H_h\mathbf H_L)
	+\beta_x\tr(\mathbf H_h\circ\mathbf H_L)
	-\beta_x\tr\{(\mathbf u_h\mathbf u_h^*)\circ\boldsymbol{\Pi}_b\}\Big],
\end{align}
which gives
\begin{equation}\label{eq:flat-bulk-k}
	\mathbf k_{n,L}=(e_\bR(h))_{h=1}^M,\qquad
	\mathbf k_{n,W}=2\delta_e^2(1+c_{nM})\mathbf k_{n,L}.
\end{equation}
The bulk variance terms reduce to
\begin{align}\label{eq:flat-bulk-variances}
	B_{n,LL}
	&=-(1+\alpha_x)\{\log(1-c_{nM})+c_{nM}\}
	+\frac1n\tr\{(1+\alpha_x)\mathbf H_L^2+\beta_x\mathbf H_L^{\circ2}\},\\
	B_{n,WW}
	&=2(1+\alpha_x)\delta_e^4c_{nM}^2
	+\frac{4\delta_e^4(1+c_{nM})^2}{n}
	\tr\{(1+\alpha_x)\mathbf H_L^2+\beta_x\mathbf H_L^{\circ2}\}.
	\label{eq:flat-bulk-variance-W}
\end{align}
Details of these kernel calculations are given in the Supplementary Material.

Logarithmic inference requires $c_n=p/n\to c\in(0,1)$, ensuring that the sample correlation matrix is asymptotically nonsingular.

\begin{corollary}[CLT for $T_L$]\label{4.1}
	Under $H_0$ in \eqref{t1} with $\mathbf G_n=\mathbf I_p$, Assumptions~\ref{ass1}--\ref{ass2}, and $c_n=p/n \rightarrow c \in (0,1)$,
	\nopagebreak[4]
	$$
	\frac{T_L-p \mathcal{I}_L -{\mu}_L}{{\sigma}_L} \xrightarrow{d} N(0,1),
	$$
	where
$
		\mathcal{I}_L=\frac{c_n-1}{c_n}\log(1-c_n)-1,$ $
		{\mu}_L=\frac{\alpha_x}{2}\log(1-c_n)+\frac{\alpha_x+1}{2}c_n,$ $
		{\sigma}^2_L=-(\alpha_x+1)\{\log(1-c_n)+c_n\}.
$
	
	Under $H_1$ in \eqref{t1} and Assumptions~\ref{ass1}--\ref{ass5}, suppose $\liminf_{n\to\infty}\hat{\sigma}_L^2>0$. Then
	$$
	\frac{s_{n,L}\left(T_L-(p-M)\hat{\mathcal{I}}_L
		-\sum_{h=1}^M\log\varphi_n(\alpha_h)\right)-\hat{\mu}_L}
	{\hat{\sigma}_L}\xrightarrow{d}N(0,1),
	$$
	where
$
		\hat{\mathcal{I}}_L
		=\frac{c_{nM}-1}{c_{nM}}\log(1-c_{nM})-1+\log\delta_e,$ $\hat{\sigma}_L^2=s_{n,L}^2\left\{
		\boldsymbol{\omega}_{n,L}^\top\mathbf{C}\boldsymbol{\omega}_{n,L}
		+2\boldsymbol{\omega}_{n,L}^\top\mathbf{k}_{n,L}+B_{n,LL}\right\},$ $
		\hat{\mu}_L=s_{n,L}\left[
		\left(\frac{\alpha_x}{2}+M\right)\log(1-c_{nM})
		-\frac{\beta_x\vartheta_{n,b}c_{nM}}{2}+b_\bR\right].
$
\end{corollary}

The quadratic statistic remains applicable throughout the regime $c_n\to c\in(0,\infty)$ of Assumption~\ref{ass2}.

\begin{corollary}[CLT for $T_W$]\label{4.2}
	Under $H_0$ in \eqref{t1} with $\mathbf G_n=\mathbf I_p$ and Assumptions~\ref{ass1}--\ref{ass2},
	\nopagebreak[4]
	$$
	\frac{T_W-p{\mathcal{I}}_w -{\mu}_w}{{\sigma}_w} \xrightarrow{d} N(0,1),
	$$
	where ${\mathcal{I}}_w=1+c_n,\  {\mu}_w=-c_n,\  {\sigma}^2_w=2(\alpha_x+1)c_n^2$.

	Under $H_1$ in \eqref{t1} and Assumptions~\ref{ass1}--\ref{ass5}, suppose $\liminf_{n\to\infty}\hat{\sigma}_w^2>0$. Then
	$$
	\frac{s_{n,W}\left(T_W-(p-M)\hat{\mathcal{I}}_w
		-\sum_{h=1}^M\varphi_n^2(\alpha_h)\right)-\hat{\mu}_w}
	{\hat{\sigma}_w}\xrightarrow{d}N(0,1),
	$$
	where
$
		\hat{\mathcal{I}}_w=\delta_e^2(1+c_{nM}),$ $
		 \hat{\sigma}_w^2=s_{n,W}^2\left\{
		\boldsymbol{\omega}_{n,W}^\top\mathbf{C}\boldsymbol{\omega}_{n,W}
		+2\boldsymbol{\omega}_{n,W}^\top\mathbf{k}_{n,W}+B_{n,WW}\right\},$ and $$
		\hat{\mu}_w=s_{n,W}\delta_e^2\left[
		(\alpha_x+\beta_x\vartheta_{n,b})c_{nM}-2(c_{nM}+2)a_\bR
		+2(2c_{nM}+3)b_\bR-Mc_{nM}(c_{nM}+2)\right].
$$
\end{corollary}

The null limits yield a level-$\epsilon$ lower-tail test based on $T_L$ and an upper-tail test based on $T_W$, with critical values $q_{n,L}=z_\epsilon\sigma_L+p\mathcal{I}_L+\mu_L$ and $q_{n,W}=z_{1-\epsilon}\sigma_w+p\mathcal{I}_w+\mu_w$. Under $H_1$, define the corresponding centers
\[
m_{n,L}=(p-M)\hat{\mathcal{I}}_L+\sum_{h=1}^M\log\varphi_n(\alpha_h),\qquad
m_{n,W}=(p-M)\hat{\mathcal{I}}_w+\sum_{h=1}^M\varphi_n^2(\alpha_h).
\]
Writing $P_L$ and $P_W$ for the rejection probabilities, the Gaussian approximations give
\begin{align*}
	P_L&\approx\Phi\!\left(\frac{s_{n,L}(q_{n,L}-m_{n,L})-\hat\mu_L}{\hat\sigma_L}\right),&
	P_W&\approx1-\Phi\!\left(\frac{s_{n,W}(q_{n,W}-m_{n,W})-\hat\mu_w}{\hat\sigma_w}\right).
\end{align*}
These expressions are used in Section~S.4.2 of the Supplementary Material to compare $T_L$, $T_W$ and Roy's largest-root test across different spike strengths.

\section{Simulation studies}\label{simula}

The simulations evaluate the two procedures developed in Section~\ref{app}. We examine the empirical size of the residual-spectrum test when factor variation is incorporated into the standard error, and compare the power of full-spectrum tests as the signal is distributed across different numbers of leading components.

\subsection{Two-sample inference for residual structure}\label{sec:reference-simulation}

The two-sample design varies residual spectral concentration while keeping total factor strength fixed. To evaluate the test in Corollary~\ref{prop:reference-bulk}, we generate two independent groups $\ell=0,1$ with $q=p-M$ according to
\begin{equation}\label{eq:reference-simulation-model}
	\begin{aligned}
		\mathbf R_\ell&=\mathbf V_\ell\operatorname{diag}
		(\alpha_{\ell,1},\ldots,\alpha_{\ell,M},t_{\ell,1},\ldots,t_{\ell,q})\mathbf V_\ell^\top,
		&\mathbf V_\ell&=\mathbf O_p\mathbf P_\ell,\\
		\mathbf Y_\ell&=\mathbf G_\ell\mathbf X_\ell,
		&\mathbf G_\ell&=\mathbf V_\ell\boldsymbol\Lambda_\ell^{1/2}\mathbf U_\ell^\top,
		\quad \mathbf U_\ell\in\{\mathbf I_p,\mathbf O_p\},
	\end{aligned}
\end{equation}
where $\mathbf O_p$ is a normalized Hadamard matrix and $\mathbf P_\ell$ is a fixed permutation. The spikes are $\alpha_{\ell,h}=(p-q\delta_e)w_{\ell,h}$, with weights summing to one. To vary residual concentration without changing the bulk mean, let $\widetilde q_n=2\lfloor q/2\rfloor$ and choose each pattern $\zeta_{\ell,j}\in\{-1,0,1\}$ with $\widetilde q_n/2$ entries of each sign and zeros elsewhere. Set
\begin{equation}\label{eq:reference-simulation-alternative}
	t_{\ell,j}=\delta_e\bigl(1+\chi_\ell p^{-1/4}\zeta_{\ell,j}\bigr),\qquad
	\Delta_{b,n}=\widetilde q_n\delta_e^2(\chi_1^2-\chi_0^2)p^{-1/2}.
\end{equation}
The primary null takes $\chi_0=\chi_1=1$, producing unequal bulk eigenvalues with equal residual concentration. Alternatives fix $\chi_0=1$ and increase $\chi_1\in\{1.25,1.5,1.75,2,2.5\}$, while the equal-bulk benchmark uses $\chi_0=\chi_1=0$. At fixed $c=p/n$, these alternatives keep the standardized signal of constant order as $p$ and $n$ increase.

We fix $M=2$, $w_{\ell,1}=5/9$ and $w_{\ell,2}=4/9$, and compare feasible and oracle standardizations. Common and permuted spectral bases isolate the effect of residual concentration from changes in eigenvectors. Gaussian and standardized $t(10)/\sqrt{5/4}$ coordinates, together with the two choices of $\mathbf U_\ell$, assess the distributional and orientation effects entering the variance.

Under equal residual concentration, the feasible test $Z_{b,n}$ has rejection rates close to the nominal 5\% level, while those of $Z'_{b,n}$ range from 25.98\% to 48.60\% (Table~\ref{tabsimulation1} and the Supplementary Material). At $p=512$, the empirical mean of $\hat V_n/V_n$ ranges from $0.9735$ to $0.9918$, showing accurate variance estimation in the larger settings.

Power increases with residual concentration and sample size (Table~\ref{tabsimulation3}). Results are similar under common and permuted spectral bases. Gaussian power is stable across coordinate orientations, while standardized $t$ coordinates with $\mathbf U=\mathbf O_p$ give power close to the Gaussian benchmark.

\begin{table}[!htbp]
	\centering
\spacingset{1}
\small
\setlength{\tabcolsep}{3.5pt}
\renewcommand{\arraystretch}{1.0}
\captionsetup{font=small,skip=5pt}
	\caption{Empirical rejection rates under $H_{0b}$ for distributions $N(0,1)$ and $t(10)/\sqrt{5/4}$ at level $\epsilon=0.05$ when $M=2$ and  $\mathbf{U}_0=\mathbf{U}_1=\mathbf{I}_p$.}
	\begin{tabular}{cc r r r r r r r r}
		\toprule
		& & \multicolumn{4}{c}{$N(0,1)$} & \multicolumn{4}{c}{$t(10)/\sqrt{5/4}$} \\
		\cmidrule(lr){3-6} \cmidrule(lr){7-10}
		&  &\multicolumn{2}{c}{$\delta_e=0.5$} & \multicolumn{2}{c}{$\delta_e=0.7$} & \multicolumn{2}{c}{$\delta_e=0.5$} & \multicolumn{2}{c}{$\delta_e=0.7$}\\
		\cmidrule(lr){3-4} \cmidrule(lr){5-6} \cmidrule(lr){7-8} \cmidrule(lr){9-10}
		$c$ & $(p,n)$ & $r_{feasible }$& $r_{bulk}$& $r_{feasible }$&$ r_{bulk}$& $r_{feasible}$& $r_{bulk}$& $r_{feasible }$&$ r_{bulk}$\\
		\midrule
		\multirow{4}{*}{0.5}
		& $(64,128)$ & 7.78 & 40.40 & 5.16 & 28.04 & 7.38 & 42.32 & 6.02 & 31.80 \\
		& $(128,256)$ & 6.68 & 42.34 & 5.86 & 33.96 & 6.36 & 43.84 & 5.98 & 36.38 \\
		& $(256,512)$ & 5.42 & 43.86 & 4.90 & 36.90 & 5.46 & 45.68 & 5.30 & 39.62 \\
		& $(512,1024)$ & 5.44 & 46.10 & 5.54 & 40.96 & 4.80 & 46.84 & 5.34 & 42.12 \\
		\multirow{4}{*}{1}
		& $(64,64)$ & 6.90 & 39.58 & 5.44 & 27.60 & 7.20 & 41.92 & 6.12 & 31.98 \\
		& $(128,128)$ & 6.30 & 41.86 & 5.90 & 32.26 & 6.54 & 43.92 & 5.48 & 36.38 \\
		& $(256,256)$ & 5.42 & 42.38 & 5.62 & 37.86 & 6.52 & 45.98 & 5.42 & 38.96 \\
		& $(512,512)$ & 5.62 & 45.84 & 5.50 & 40.58 & 5.58 & 46.76 & 5.18 & 42.66 \\
		\multirow{4}{*}{2}
		& $(128,64)$ & 7.08 & 42.04 & 5.86 & 31.06 & 6.20 & 41.58 & 6.16 & 35.42 \\
		& $(256,128)$ & 6.10 & 43.82 & 5.52 & 37.20 & 6.30 & 45.46 & 5.82 & 39.22 \\
		& $(512,256)$ & 5.26 & 46.38 & 5.98 & 42.40 & 5.90 & 48.26 & 5.54 & 42.78 \\
		& $(1024,512)$ & 5.44 & 46.64 & 5.34 & 43.98 & 5.76 & 48.60 & 5.48 & 44.10 \\
		\bottomrule
	\end{tabular}
	\label{tabsimulation1}
\end{table}

\begin{table}[!htbp]
\centering
\spacingset{1}
\small
\setlength{\tabcolsep}{3.5pt}
\renewcommand{\arraystretch}{1.0}
\captionsetup{font=small,skip=5pt}
\caption{Empirical powers for testing $H_{0b}$ at level $\epsilon=0.05$, with $\chi_0=1$, increasing $\chi_1$, $M=2$, $\mathbf{U}_0=\mathbf{U}_1=\mathbf{I}_p$ and common spectral bases.}
\label{tabsimulation3}
\begin{tabular*}{\textwidth}{@{\extracolsep{\fill}}cc*{10}{r}@{}}
\toprule
& & \multicolumn{5}{c}{$N(0,1)$} & \multicolumn{5}{c}{$t(10)/\sqrt{5/4}$} \\
\cmidrule(lr){3-7}\cmidrule(lr){8-12}
& & \multicolumn{5}{c}{$\chi_1$} & \multicolumn{5}{c}{$\chi_1$} \\
$c$ & $(p,n)$ & 1.25 & 1.5 & 1.75 & 2 & 2.5 & 1.25 & 1.5 & 1.75 & 2 & 2.5 \\
\midrule
\multirow{4}{*}{0.25} & $(64,256)$ & 17.52 & 38.86 & 67.08 & 88.84 & 99.78 & 15.82 & 33.02 & 56.56 & 79.00 & 98.50 \\
 & $(128,512)$ & 17.80 & 40.40 & 69.80 & 91.34 & 99.98 & 16.00 & 31.72 & 58.14 & 81.16 & 99.00 \\
 & $(256,1024)$ & 16.94 & 39.76 & 72.00 & 92.92 & 100.00 & 14.72 & 32.14 & 58.72 & 83.70 & 99.48 \\
 & $(512,2048)$ & 17.28 & 41.32 & 73.26 & 94.40 & 99.94 & 13.96 & 31.68 & 59.22 & 83.92 & 99.72 \\
\addlinespace[3pt]
\multirow{4}{*}{0.5} & $(64,128)$ & 12.64 & 21.84 & 40.56 & 59.06 & 89.40 & 12.76 & 20.56 & 33.18 & 48.64 & 79.84 \\
 & $(128,256)$ & 11.70 & 22.78 & 38.54 & 61.10 & 92.00 & 11.32 & 18.48 & 30.62 & 49.16 & 83.36 \\
 & $(256,512)$ & 11.22 & 21.28 & 39.08 & 62.00 & 93.56 & 9.92 & 18.90 & 31.22 & 48.16 & 84.30 \\
 & $(512,1024)$ & 10.90 & 21.22 & 39.14 & 62.44 & 94.98 & 9.72 & 17.42 & 30.96 & 50.20 & 85.40 \\
\addlinespace[3pt]
\multirow{4}{*}{1} & $(64,64)$ & 10.04 & 14.18 & 21.42 & 29.92 & 54.02 & 11.04 & 13.94 & 19.40 & 28.06 & 46.96 \\
 & $(128,128)$ & 8.96 & 14.04 & 19.36 & 29.64 & 55.00 & 8.28 & 13.38 & 18.74 & 25.88 & 45.66 \\
 & $(256,256)$ & 7.84 & 12.54 & 20.06 & 28.54 & 57.88 & 8.22 & 11.28 & 17.12 & 25.58 & 46.40 \\
 & $(512,512)$ & 7.74 & 11.86 & 19.28 & 30.26 & 55.74 & 8.04 & 12.24 & 17.48 & 23.48 & 45.86 \\
\bottomrule
\end{tabular*}
\end{table}

\subsection{One-sample tests of an identity correlation matrix}\label{sec:classical-numerics}

We compare the logarithmic and quadratic tests with Roy's largest-root test $T_R$ to examine how the allocation of factor signal affects power. Critical values for $T_R$ are obtained from the Tracy--Widom approximation described in the Supplementary Material.

Alternatives to $H_0:\bR=\bI_p$ use $K\in\{1,2,5\}$ spikes:
\[
\bR^{(K)}=\mathbf O_p\operatorname{diag}
\left(\alpha_1,\ldots,\alpha_K,
\underbrace{\delta_e,\ldots,\delta_e}_{p-K}\right)\mathbf O_p^*,
\qquad
\sum_{h=1}^K\alpha_h+(p-K)\delta_e=p,
\]
where $\mathbf O_p$ is a normalized $p\times p$ Hadamard matrix. Models 1.1, 1.2 and 1.3 correspond to $K=1,2,5$. For $K=2$, we take $\alpha_2/\alpha_1\in\{0.75,0.8,0.85,0.9\}$, while for $K=5$, $(\alpha_1,\ldots,\alpha_5)=\alpha_1(1,0.95,0.9,0.85,0.8)$. The trace constraint determines $\alpha_1$ for each $\delta_e$.

We draw i.i.d. $x_{ij}$ from $D_1=N(0,1)$, $D_2=t(10)/\sqrt{5/4}$ or $D_3=U(-\sqrt{3},\sqrt{3})$. With $p/n=0.5$ and $p=64,128,256$, each configuration uses 5000 replications. Size is evaluated at $\epsilon=0.05$ and $0.01$, and power at $\epsilon=0.05$ for $T_L$, $T_W$ and $T_R$ under Models 1.1, 1.2 and 1.3.

Across all three distributions, the null rejection rates of $T_L$ and $T_W$ are close to the nominal levels (Table~\ref{testH0}), supporting the power comparisons below.

\begin{table}[!htbp]
	\centering
\spacingset{1}
\small
\setlength{\tabcolsep}{3.5pt}
\renewcommand{\arraystretch}{1.0}
\captionsetup{font=small,skip=5pt}
	\caption{Empirical rejection rates under $H_0$ for the three coordinate distributions at levels $\epsilon=0.05$ and $0.01$.}
	\label{testH0}
	\begin{tabular*}{\textwidth}{@{\extracolsep{\fill}}cc*{6}{c}@{}}
		\toprule
		& & \multicolumn{3}{c}{$\epsilon=0.05$} & \multicolumn{3}{c}{$\epsilon=0.01$} \\
		\cmidrule(lr){3-5}\cmidrule(lr){6-8}
		test & $(p,n)$ & $D_1$ & $D_2$ & $D_3$ & $D_1$ & $D_2$ & $D_3$ \\
		\midrule
		\multirow{3}[1]{*}{$T_L$} & (64,128) & 4.82 & 5.28 & 5.26 & 0.74 & 1.10 & 1.18 \\
		& (128,256) & 5.00 & 5.08 & 5.38 & 1.30 & 0.80 & 1.14 \\
		& (256,512) & 5.48 & 4.58 & 5.12 & 1.20 & 1.06 & 1.04 \\
		\multirow{3}[0]{*}{$T_W$} & (64,128) & 4.52 & 4.94 & 4.48 & 1.10 & 1.40 & 1.06 \\
		& (128,256) & 5.16 & 4.80 & 5.52 & 1.10 & 1.18 & 1.04 \\
		& (256,512) & 4.82 & 4.88 & 4.96 & 1.10 & 0.96 & 1.18 \\
		\multirow{3}[1]{*}{$T_R$} & (64,128) & 0.98 & 0.84 & 0.84 & 0.12 & 0.10 & 0.14 \\
		& (128,256) & 1.40 & 1.46 & 1.66 & 0.22 & 0.24 & 0.18 \\
		& (256,512) & 1.82 & 1.90 & 1.70 & 0.32 & 0.28 & 0.30 \\
		\bottomrule
	\end{tabular*}
\end{table}

With a single dominant spike, $T_R$ has the highest empirical power. For two spikes, $T_W$ and $T_R$ outperform $T_L$, with little variation as $\alpha_2/\alpha_1$ ranges from $0.75$ to $0.9$. Power increases with dimension and spike strength, with similar patterns across the three coordinate distributions. Detailed results for one and two spikes are reported in the Supplementary Material.

With five comparable spikes, $T_W$ has the highest power (Table~\ref{m5}). At $(p,n)=(64,128)$ and $\delta_e=0.92$, its power ranges from 98.92\% to 99.14\%, compared with 89.22\% to 89.66\% for $T_R$. These results show that aggregating signal across several leading components can substantially improve detection when the factor signal is distributed across multiple directions.

\begin{table}[!htbp]
\centering
\spacingset{1}
\small
\setlength{\tabcolsep}{3.5pt}
\renewcommand{\arraystretch}{1.0}
\captionsetup{font=small,skip=5pt}
	\caption{Empirical powers for Model 1.3 under the three coordinate distributions at level $\epsilon=0.05$. Within each distribution, the columns correspond to $\delta_e=0.92,0.90,0.88$.}
	\label{m5}
	\begin{tabular*}{\textwidth}{@{\extracolsep{\fill}}cc*{9}{c}@{}}
		\toprule
		& & \multicolumn{3}{c}{$D_1$} & \multicolumn{3}{c}{$D_2$} & \multicolumn{3}{c}{$D_3$} \\
		\cmidrule(lr){3-5}\cmidrule(lr){6-8}\cmidrule(lr){9-11}
		test & $(p,n)$ & .92 & .90 & .88 & .92 & .90 & .88 & .92 & .90 & .88 \\
		\midrule
		\multirow{3}[1]{*}{$T_L$} & (64,128) & 80.64 & 96.62 & 99.86 & 81.06 & 97.20 & 99.90 & 80.00 & 96.74 & 99.86 \\
		& (128,256) & 100 & 100 & 100 & 100 & 100 & 100 & 100 & 100 & 100 \\
		& (256,512) & 100 & 100 & 100 & 100 & 100 & 100 & 100 & 100 & 100 \\
		\multirow{3}[0]{*}{$T_W$} & (64,128) & 98.92 & 100 & 100 & 99.14 & 100 & 100 & 99.08 & 100 & 100 \\
		& (128,256) & 100 & 100 & 100 & 100 & 100 & 100 & 100 & 100 & 100 \\
		& (256,512) & 100 & 100 & 100 & 100 & 100 & 100 & 100 & 100 & 100 \\
		\multirow{3}[1]{*}{$T_R$} & (64,128) & 89.62 & 99.50 & 99.98 & 89.66 & 99.54 & 100 & 89.22 & 99.62 & 99.98 \\
		& (128,256) & 100 & 100 & 100 & 100 & 100 & 100 & 100 & 100 & 100 \\
		& (256,512) & 100 & 100 & 100 & 100 & 100 & 100 & 100 & 100 & 100 \\
		\bottomrule
	\end{tabular*}
\end{table}

\section{Residual genomic dependence in HapMap data}\label{sec:empirical}

We apply the two-sample procedure to compare residual genomic dependence between male and female groups after accounting for prominent ancestry components. The analysis uses autosomal genotypes from the 1000-SNP HapMap subset included in \href{https://bioconductor.org/packages/3.22/bioc/src/contrib/SNPRelate_1.44.0.tar.gz}{\texttt{SNPRelate} 1.44.0}, with population and sex labels from the \href{https://ftp.1000genomes.ebi.ac.uk/vol1/ftp/release/20130502/integrated_call_samples_v3.20130502.ALL.panel}{1000 Genomes Phase 3 sample panel}.  The two groups have identical ancestry composition, each with $n=86$ individuals from the YRI, CEU, CHB and JPT populations. We retain $p=110$ complete autosomal markers with common allele frequencies across populations, selected to capture ancestry variation while limiting strong local dependence.

Both groups exhibit two clearly separated leading eigenvalues. Writing the sample correlation eigenvalues as $\widehat\lambda_{\ell,1}\geq\cdots\geq\widehat\lambda_{\ell,p}$, the first two are $(26.57,12.24)$ for males and $(24.67,11.94)$ for females, with second-to-third eigenvalue ratios of $4.25$ and $4.09$. We therefore take $M=2$ and $q=108$. To assess the equal-eigenvalue reference for the residual spectrum, we invert the Mar\v{c}enko--Pastur equation using the residual empirical distribution. Across the three imaginary parts considered, the recovered spectrum has relative mean absolute deviations between $0.009$ and $0.027$, with at least $98.5\%$ of its mass lying within $10\%$ of its mean. Gaussian benchmarks at the same $(n,q)$ show comparable behavior, supporting the equal-bulk approximation for the residual comparison.

With males labeled as group $0$ and females as group $1$, define the residual concentration $Q_\ell=\sum_{j=3}^{p}\widehat\lambda_{\ell,j}^2$ and the between-group contrast $D_b=Q_1-Q_0$. Corollary~\ref{prop:reference-bulk} gives the feasible variance estimate
\[
\widehat V_\ell=4\widehat\delta_\ell^4(1+q/n)^2
\frac{\widehat\nu_\ell}{n},\qquad
\widehat\nu_\ell=\frac1n\sum_{r=1}^n(e_{\ell r}-p)^2,
\]
where $e_{\ell r}$ is the squared length of an observation after coordinatewise sample standardization.

For comparison, we use the same numerator with a separate-bulk variance
\[
\widehat V_{\mathrm{sep}}
=8\bar\delta^4(q/n)^2,
\qquad
\bar\delta=\frac{1}{2q}\sum_{\ell=0}^{1}
\sum_{j=3}^{p}\widehat\lambda_{\ell,j}.
\]
This benchmark treats the residual eigenvalues scaled by $\bar\delta$ as a sample correlation spectrum with identity population correlation, while $\widehat V_\ell$ retains the contribution of factor variation to the sample-variance normalization.

\begin{table}[!htbp]
	\centering
	\spacingset{1}
	\small
	\setlength{\tabcolsep}{3.5pt}
	\renewcommand{\arraystretch}{1.0}
	\captionsetup{font=small,skip=5pt}
	\caption{Residual quadratic sums and two standardizations of $D_b=6.17$.}\label{tab:hapmap-empirical}
	\begin{tabular}{lrrrr}
		\toprule
		Quantity or method & Male & Female & Combined $Z$ & Two-sided $p$ \\
		\midrule
		$Q_\ell$ & 108.03 & 114.20 & & \\
		$\widehat V_\ell$ & 21.50 & 20.42 & & \\
		Proposed variance & & & 0.953 & 0.341 \\
		Separate correlation LSS variance & & & 3.876 & $1.06\times10^{-4}$ \\
		\bottomrule
	\end{tabular}
\end{table}

Accounting for factor variation gives a standard error of $6.47$ for the observed difference $D_b=6.17$, compared with $1.59$ under the separate correlation LSS benchmark (Table~\ref{tab:hapmap-empirical}). The proposed test yields $Z=0.953$ and $p=0.341$, consistent with equal residual concentration at the 5\% level, while the benchmark rejects. Across selections of the 80, 90, 100 and 110 most informative markers, the proposed two-sided $p$-values range from $0.245$ to $0.341$, compared with $2.24\times10^{-6}$ to $1.06\times10^{-4}$ for the benchmark. The difference shows that ancestry-related factor variation materially affects uncertainty in the residual comparison.

\section*{Disclosure Statement}
The authors report there are no competing interests to declare.

\section*{Data Availability Statement}
The empirical analysis uses the publicly available 1000-SNP HapMap subset included in \texttt{SNPRelate} version 1.44.0, with population and sex labels obtained from the 1000 Genomes Phase 3 sample panel. Data sources is described in Section~\ref{sec:empirical}.

\end{document}